\documentclass[12pt]{article}
\usepackage{amsfonts}
\usepackage{amsthm}
\usepackage{graphicx}
\usepackage{latexsym}
\usepackage{amssymb}
\usepackage{amsmath}
\usepackage{colortbl}
\definecolor{text1}{cmyk}{1,.65,0,0} 
\definecolor{text2}{rgb}{1,0,0} 
\definecolor{text3}{cmyk}{0,0,0,1} 
\definecolor{text4}{cmyk}{0,0,0,0.5} 
\newcommand{\blue}[1]{{\color{text1}#1}}

\newtheorem{definition}{Definition}[section]
\newtheorem{theorem}{Theorem}[section]
\newtheorem{lemma}{Lemma}[section]
\newtheorem{corollary}{Corollary}[section]

\makeatletter

\begin{document}

\begin{center}
   {\bf On efficient prediction and predictive density estimation for spherically symmetric models  \footnote{\today}}
\end{center}

\begin{center}
{\sc Dominique Fourdrinier$^{a}$, \'Eric Marchand$^{b}$, William E. Strawderman$^{c}$} \\

{\it a  Universit\'{e} de Normandie, INSA Rouen, UNIROUEN, UNIHAVRE, LITIS, avenue de l'Universit\'{e}, BP 12, 76801 Saint-\'{E}tienne-du-Rouvray,  FRANCE  \quad (e-mail: dominique.fourdrinier@univ-rouen.fr) } \\
    
{\it b  Universit\'e de
    Sherbrooke, D\'epartement de math\'ematiques, Sherbrooke Qc,
    CANADA, J1K 2R1 \quad (e-mail: eric.marchand@usherbrooke.ca) } \\

{\it c  Rutgers University, Department of Statistics and Biostatistics, 501 Hill
Center, Busch Campus, Piscataway, N.J., USA, 08855 \quad (e-mail:
straw@stat.rutgers.edu) }
\end{center}
\vspace*{0.2cm}
\begin{center}
{\sc Summary} \\
\end{center}
\vspace*{0.2cm}
\small
	Let $X,U,Y$ be spherically symmetric distributed having density $$\eta^{d +k/2} \, f\left(\eta(\|x-\theta|^2+ \|u\|^2  + \|y-c\theta\|^2 ) \right)\,,$$ with unknown parameters $\theta \in \mathbb{R}^d$ and  $\eta>0$, and with known density $f$ and constant $c >0$. Based on observing $X=x,U=u$, we consider the problem of obtaining a predictive density $\hat{q}(y;x,u)$ for $Y$ as measured by the expected Kullback-Leibler loss.
A benchmark procedure is the minimum risk equivariant density $\hat{q}_{mre}$, which is Generalized Bayes with respect to the prior $\pi(\theta, \eta) = \eta^{-1}$.     For $d \geq 3$, we obtain improvements on  $\hat{q}_{mre}$, and further show that the dominance holds simultaneously for all $f$ subject to finite moments and finite risk conditions.  We also obtain that the Bayes predictive density with respect to the harmonic prior $\pi_h(\theta, \eta) =\eta^{-1} \|\theta\|^{2-d}$ dominates  $\hat{q}_{mre}$ simultaneously for all scale mixture of normals $f$.   \\

\noindent  The results hinges on duality with a point prediction problem, as well as posterior representations for $(\theta, \eta)$, which are of interest on their own.  Namely, we obtain for 
$d \geq 3$, point predictors $\delta(X,U)$ of $Y$ that dominate the benchmark predictor $cX$ simultaneously for all $f$, and simultaneously for risk functions $\mathbb{E}_f \left[\rho \,( \|Y - \delta(X,U) \|^2  + (1+c^2)\|U\|^2)\right]$, with $\rho$ increasing and concave on $\mathbb{R}_+$, and including the squared error case $\mathbb{E}_f \left[\,( \|Y - \delta(X,U) \|^2  \right]$    

\vspace*{0.5cm}
\small
\noindent  {\it AMS 2010 subject classifications:}   62C20, 62C86, 62F10, 62F15, 62F30

\noindent {\it Keywords and phrases}:  Bayes estimators; Dominance; Duality; Kullback-Leibler; Multivariate normal; Multivariate Student; Plug-in; Point prediction; Predictive densities; Scale mixture of normals; Spherically symmetric.

\normalsize
\section{Introduction and preliminary results}

\noindent {\bf A.} Prediction is at the heart of statistics, but the study of the efficiency of prediction methods often takes a back seat to estimation.  There is perhaps a reason for this.  Indeed, consider $Z_1,Z_2 \in \mathbb{R}^d$ independently and identically distributed (i.i.d.) random variables, with $\mathbb{E}(Z_1)=\theta$, $\hbox{Cov}(Z_1)= \Sigma$ , and the problem of predicting $Z_2$ based on $Z_1$.  If our prediction is $\delta(Z_1)$ and the penalty is squared error, then
$$  \mathbb{E}[\|Z_2-\delta(Z_1)\|^2]\,=\, tr \Sigma \, + \, \mathbb{E}[\|\delta(Z_1)-\theta\|^2]\,, \hbox{for all } \theta, $$
so that the frequentist squared error risk of $\delta(Z_1)$ as a predictor of $Z_2$ is determined by its frequentist risk as a point estimator of $\theta$.  For instance, in the case of the distribution of $Z_1,Z_2$ being a multivariate normal distribution with $d\geq 3$, shrinkage or Stein-type estimators $\delta(Z)$ that dominate $Z_1$ as estimators of $\theta$ (e.g., Strawderman 2003) yield improved predictors $\delta(Z_1)$ of $Z_2$ as described above. \footnote{Similarly, if the penalty is given by $(Z_2-\delta(Z_1))' Q (Z_2-\delta(Z_1))$ with $Q$ positive definite, then
we have the decomposition 
$$ tr Q \Sigma \, + \, \mathbb{E} \left((\delta(Z_1) - \theta)'Q (\delta(Z_1) - \theta) \right)\,$$
and another clear correspondence with a familiar point estimation problem. }
If the prediction penalty is not squared error, then the above correspondence is obviously different and relationships between prediction and estimation are more subtle.   Moreover, the decision-maker may well wish to select an alternative to squared error penalty and, namely, a penalty that is non-convex or bounded, or both.   \\

\noindent  {\bf B.}  Alternatively, predictive density estimation aims at providing 
the richest description of an unobserved random variable in the form of a predictive density over the domain of possible values.   One obtains a surrogate density for a future or missing value, based on current or historical data.    Bayesian strategies for deriving predictive densities can be naturally formulated in response to a given prior and a measure of divergence between densities, such as Kullback-Leibler.  There also arise issues of efficiency and frequentist risk evaluation of predictive densities.
In this regard, following seminal contributions such as Aitchison (1975), Aitchison and Dunsmore (1975), and Komaki (2001), further challenges relative to the efficiency of predictive densities, for various models and loss functions, have generated much more recent interest, as exemplified by the work of Liang and Barron (2004), George, Liang and Xu (2006), Komaki (2006, 2007), Brown, George and Xu (2008), Fourdrinier et al. (2011), and many others including those referred to below.   Namely, several parallels between point and predictive density estimation have surfaced (e.g., the inadmissibility of the minimum risk equivariant procedures for squared error and Kullback-Leibler losses, for normal observables in three dimensions or more), including Bayesian procedures.  However, this is less the case for connections between point prediction and predictive density estimation.   
As well, for general spherically symmetric models with unknown location and scale parameters, including the normal model, much less is known on the efficiency of predictive density estimators, Bayesian or otherwise.   \\

\noindent {\bf C.}  This paper's contributions relate to both point prediction and predictive density estimation, as well as connections which generate further findings for the latter.  We consider broadly a predictive density estimation problem based on

\begin{equation}
\label{modelXUY}
X,U,Y|\theta,\eta \sim \eta^{d +k/2} \, f\left(\eta(\|x-\theta|^2+ \|u\|^2  + \|y-c\theta\|^2 ) \right)\,,
\end{equation}
with $x,y, \theta \in \mathbb{R}^{d}$, $u \in \mathbb{R}^k$, $\eta^{-1/2}$ a scale parameter,  $c$ positive and known, and $f$ a known spherically symmetric density on $\mathbb{R}^{2d+k}$.   
Such a model arises quite generally as a canonical form generated from a linear model (see for instance Fourdrinier and Strawderman, 2010).  It includes the multivariate normal model with independent components $X$, $U$, and $Y$ 

\begin{equation}\label{modelXYS}
X \sim N_d(\theta,  \eta^{-1} I_{d}), \;Y \sim N_d(c\theta,  \eta^{-1} I_{d}), \; S=U'U \sim \eta^{-1} \chi^{2}_{k}, \;\mbox{independent},
\end{equation}
where the objective is to predict $Y$ based on $(X,S)$.  Model (\ref{modelXYS}) applies for the familiar set-up where we observe  $X_1, \ldots, X_n$ independently distributed $N_{d}(\mu, \sigma^2 I_{d})$ and wish to predict $Y$ as above.  This is achieved by setting $X=\sqrt{n} \bar{X}$, $S=\sum_{i=1}^n \|X_i-\bar{X}\|^2$, $\theta=\sqrt{n} \mu$, $c=n^{-1/2}$, and $k=(n-1){d}$.  Otherwise, model (\ref{modelXUY}) encapsulates situations where the signals $X$, $Y$ are not independent of the residual vector $U$ and exhibit a spherically symmetric dependence.  \\

\noindent  Based on $(X,U)$, we seek efficient predictive densities $\hat{q}(y;x,u), y \in \mathbb{R}^{d}$, for the conditional density $q_{\theta, \eta}(\cdot|x,u)$ of $Y$ given $x,u$.  We evaluate the performance of such predictive densities by Kullback-Leibler loss
\begin{equation}
\label{Kullback-Leibler}
L_{KL}((\theta, \eta), \hat{q})  = \int_{\mathbb{R}^{d}}  q_{\theta, \eta}(y|x,u) \, \log \left(\frac{q_{\theta, \eta}(y|x,u)}{\hat{q}(y; x,u)}  \right)  \, dy \,,
\end{equation}
and associated frequentist risk taken with respect to the marginal density $p_{\theta, \eta}$ of $X,U$, given by
\begin{eqnarray}
\nonumber
R_{KL}((\theta, \eta), \hat{q})\, & = & \int_{\mathbb{R}^{d+k}} L_{KL}((\theta, \eta), \hat{q}) \, p_{\theta, \eta}(x,u) \, dx \, du \\
  & = &
\label{riskKL}  \mathbb{E}^{X,U,Y} \, \log \left(\frac{ \, q_{\theta, \eta}(Y|X,U)}{\hat{q}(Y; X,U)}  \right)
\,   .
\end{eqnarray}

\noindent {\bf D.}
A benchmark predictive density is the Bayes predictive density estimator $\hat{q}_{\pi_0,}(\cdot;X,U)$ with respect to the prior measure $\pi_0(\theta, \eta) = \frac{1}{\eta}$.  It is also minimax and the minimum risk equivariant (mre) predictive density with respect to changes of location and scale (e.g., Kubokawa et al., 2013).   It will be shown that it is given by a multivariate Student density, that is 

\begin{equation}
\label{mre}
\hat{q}_{\pi_0}(\cdot;(x,u)) \sim T_{d}(k,cx, \sqrt{\frac{((1+c^2) \|u\|^2}{k}})\,.
\end{equation}
 Hereafter, we refer to multivariate Student densities as follows.
\begin{definition}
\label{student}
 A $d$-variate Student distribution with degrees of freedom $\nu$, location parameter $\xi$, scale parameter $\sigma$, denoted $T_d(\nu,\xi,\sigma)$ has density given by

\begin{equation}
\label{studentpdf}
\frac{1}{\sigma^d} \, \frac{\Gamma(\frac{\nu+d}{2})}{\Gamma(\frac{\nu}{2}) (\pi \nu )^{d/2}} \, \left( 1 + \frac{\|t-\xi\|^2}{\nu \sigma^2}  \right)^{-\frac{d+\nu}{2}}\,, \, t \in \mathbb{R}^d.
\end{equation} 
\end{definition}

\noindent   For the normal case as  in (\ref{modelXYS}), the predictive density $\hat{q}_{\pi_0}$ was obtained in Aitchison and Dunsmore (1975), and shown to be minimax by Liang and Barron (2004).   However, the Bayes predictive density $\hat{q}_{\pi_0}$ is known to be inadmissible for $d \geq 3$ in the normal case.  Indeed, Kato (2009) showed that it was uniformly improved with respect to Kullback-Leibler risk by the Bayes predictive density estimator associated with the harmonic prior $\pi_h(\theta, \eta) = \eta^{-1} \|\theta\|^{2-d}$.   Moreover, further improvements (still in the normal case), even for some cases with $d<3$, were obtained by Boisbunon and Maruyama (2014), and earlier work by Komaki (2006, 2007) established the inadmissibility of $\hat{q}_{\pi_0}$ in an asymptotic framework. 

\noindent {\bf E.} Expression (\ref{mre}) is proved in {\bf Section 3.2}.  Moreover, we point out that the predictive density $\hat{q}_{\pi_0}$ does not depend on the model density $f$ in (\ref{modelXUY}) (and consequently matches the normal case solution).  In {\bf Section 3}, we elaborate on this phenomenon from a more general perspective, where a class of Bayesian inference methods, associated with separable priors of the form $\pi_1(\theta) \, \eta^a$, do not depend on $f$; and dominance results that hold in the normal case carry-over to the whole class of scale mixture of normals.  \\

\noindent  A main focus of this paper is on providing improvements on $\hat{q}_{\pi_0}$ applicable to model densities in (\ref{modelXUY}).  Bayesian solutions are presented in {\bf Section 4}.  Namely, we prove that, for $d \geq 3$ and a given scale mixture of normals $f$ in (\ref{modelXUY}), the Bayesian predictive density $\hat{q}_{\pi_h}$ with respect to the harmonic prior $\pi_h$ dominates the mre predictive density $\hat{q}_{\pi_0}$.  Moreover, both $\hat{q}_{\pi_h}$ and $\hat{q}_{\pi_0}$ do not vary with the scale mixture and the dominance holds simultaneously for all scale mixtures.  \\

\noindent As presented in {\bf Section 2}, our findings include dominating predictive densities for $d \geq 3$, which are multivariate Student densities of the form $T_d(k, c \hat{\theta}(x,u), \sqrt{\frac{(1+c^2) \|u\|^2}{k}}) $, and where $\hat{\theta}(x,u)$ is a point estimator of $\theta$.  The focus on such predictive densities, which we find convenient to denote by $q_{\pi_0,\hat{\theta}}$, leads to a key duality result presented in Lemma \ref{duality}.
More precisely, the Kullback-Leibler risk performance of predictive density $q_{\pi_0,\hat{\theta}}(\cdot; X,U)$ hinges on the performance of $c \hat{\theta}(X,U)$ as a point predictor of $Y$ under ``loss''
\begin{equation}
\label{lossrho}
\rho \left( \|Y-c\hat{\theta}(X,U)\|^2 + (1+c^2) \|U\|^2 \right)  \,,
\end{equation}
with $\rho(t) = \log(t)\,,t >0$.\footnote{While it is standard to require that a loss function be bounded below, we will nonetheless refer to $\rho(t)=\log(t)$ in (\ref{lossrho}) as a loss even though it is not so bounded.}  A general dominance result for the point prediction problem, applicable for $\rho$ increasing and concave and $d \geq 3$,  is obtained with Theorem \ref{main} and leads immediately to the predictive density estimation finding of Theorem \ref{dominancepde}.   Hence, the findings of Section 2 are contributions to both: {\bf (A)}  point prediction of $Y$ for model (\ref{modelXUY}) and losses (\ref{lossrho}), as well as {\bf (B)}  predictive density estimation of conditional density $q_{\theta, \eta}(\cdot|x,u)$ of $Y$ under Kullback-Leibler loss.  Moreover, the dominance results for {\bf (A)} are shown to be, subject to risk-finiteness, robust with respect to model density $f$ in 
{\bf (1,2)} and loss (\ref{lossrho}) for general $\rho$, while those for {\bf (B)} are shown to be robust with respect to $f$.   The techniques used to derive the classes of dominating procedures involve a Stein identity for spherical densities (Lemma \ref{Steingeneral}), as well as a concave inequality technique analogous to earlier point estimation work of Brandwein and Strawderman (1980), Brandwein and Strawderman (1991), 
Brandwein, Ralescu and Strawderman (1993), and Kubokawa, Marchand and Strawderman (2015).

\section{Results for point prediction with predictive density estimation implications}

\noindent We begin by connecting the predictive density estimation problem with a point prediction problem.

\begin{lemma}
\label{duality}
For a spherically symmetric model as in (\ref{modelXUY}),
the predictive density $q_{\pi_0,\hat{\theta}} \sim T_{d}(k, c \hat{\theta}(X,U), \sqrt{\frac{(1+c^2) \|U|^2}{k}})$  dominates the predictive density $\hat{q}_{\pi_0}$ given in (\ref{mre}) under Kullback-Leibler loss if and only if

\begin{equation}
\label{conditionlog}
\mathbb{E}_f \left[ \log \left( \|Y-cX\|^2 + (1+c^2) \|U\|^2 \right)  \right]
\geq \mathbb{E}_f \left[ \log ( \|Y-c\hat{\theta}(X,U)\|^2 + (1+c^2) \|U\|^2)  \right]\,,
\end{equation}
for all $\theta, \eta$, with strict inequality for at least one $(\theta,\eta)$.
\end{lemma}
\noindent {\bf Proof.}
We have as a difference in risks
\begin{eqnarray*}
R_{KL}((\theta,\eta), \hat{q}_{\pi_0}) - R_{KL}((\theta,\eta), q_{\pi_0, \hat{\theta}}) & = &  
\mathbb{E}_f \, log \left(\frac{q_{\theta, \eta}(Y|X,U)}{\hat{q}_{\pi_0,X}(Y)}   \right) - \mathbb{E}_f \, log \left(\frac{q_{\theta, \eta}(Y|X,U)}{\hat{q}_{\pi_0,\hat{\theta}(X,U)}(Y)} \right)
\\ 
& = & \mathbb{E}_f \, \log \left( \frac{q_{\pi_0,\hat{\theta} (X,U)}(Y)}{\hat{q}_{\pi_0,X}(Y)} \right) \\
\, & = & \mathbb{E}_f \, \log \frac{\left( 1 + \frac{\|Y-c \hat{\theta}(X,U)\|^2}{(1+c^2) \|U\|^2}  \right)^{-\frac{d+k}{2}}}{\left( 1 + \frac{\|Y-c X\|^2}{(1+c^2) \|U\|^2}  \right)^{-\frac{d+k}{2}}}  \\
\, & = & \frac{d+k}{2}  \left( \mathbb{E}_f \left[ \log \left( \|Y-cX\|^2 + (1+c^2) \|U\|^2 \right)  \right] \right)  \\
\, & - &  \frac{d+k}{2} \left( \mathbb{E}_f \left[ \log \left( \|Y-c\hat{\theta}(X,U)\|^2 + (1+c^2) \|U\|^2 \right)  \right] \right)\,,
\end{eqnarray*} 
which establishes the result.  \qed

In the following, for a vector valued function $g(t_1,t_2)$ with $\hbox{dim } g(t_1,t_2)= \hbox{dim } t_1$,
$\hbox{div}_{t_1} \, g(t_1,t_2)$ represents the divergence with respect to $t_1 $.   We will make use of the following useful identity, a version of which can be found in Fourdrinier, Strawderman and Wells (2003), obtained by integration by parts and reducing to the celebrated Stein identity in the normal case.

\begin{lemma}
\label{Steingeneral}
Let $Z \in \mathbb{R}^d, U \in \mathbb{R}^k$ have joint density $f(\|z\|^2+ \|u\|^2)$ and let $w \in \mathbb{R}^d$ be fixed.  Set $F(t)=\frac{1}{2} \int_{t}^{\infty} f(u) \, du$, and $\gamma=\int_{\mathbb{R}^{d+k}} F(\|z\|^2+ \|u\|^2)\, du \,dz$ assuming the integral is finite.  Then, we have 
for weakly differentiable $g: \mathbb{R}^{2d+k} \to \mathbb{R}^d$ and $h: \mathbb{R}^{2d+k} \to \mathbb{R}^k$
\begin{eqnarray*}
\int_{\mathbb{R}^{d+k}} z' g(z,u,w) \,  f(\|z\|^2+ \|u\|^2)\, du \,dz &=& \frac{1}{\gamma}\int_{\mathbb{R}^{p+k}} \hbox{div}_z g(z,u,w) \, F(\|z\|^2+ \|u\|^2)\, du \,dz   \\
\int_{\mathbb{R}^{d+k}} u' h(z,u,w) \,  f(\|z\|^2+ \|u\|^2)\, du \,dz &=& \frac{1}{\gamma} \int_{\mathbb{R}^{p+k}} \hbox{div}_u \, h(z,u,w) \, F(\|z\|^2+ \|u\|^2)\, du \,dz \,,
\end{eqnarray*}
provided the integrals exist.
\end{lemma}

	We now are ready to present, establish, and comment on the main point prediction result, which follows.

\begin{theorem}
\label{main}
Let $\tilde{X},\tilde{Y}  \in \mathbb{R}^{d}$, $\tilde{U} \in \mathbb{R}^k$ have joint density proportional to
\begin{equation}
\label{modelX'Y'U'}
 \eta_1^{d+k/2} f \left( \eta_1 \left( \|\tilde{x}-\mu\|^2+ \frac{\|\tilde{y}-\mu\|^2}{\beta} + 
 \frac{\|\tilde{u}\|^2}{1+\beta} \right) \right) \,,
 \end{equation} 
where $\beta>0$ (known) and $d \geq 3$.  Consider predicting $\tilde{Y}$ with $\delta(\tilde{X},\tilde{U})$ under loss $\rho(\|\delta-\tilde{Y}\|^2 + \|\tilde{U}\|^2)$, where $\rho(\cdot)$ is absolutely continuous,  increasing and concave.   Then, the predictor $\tilde{X} + \frac{\alpha \|\tilde{U}\|^2}{k+2} g(\tilde{X})$ dominates $\tilde{X}$  provided $\mathbb{E}_{\theta, \eta_1} \|g(\tilde{X})\|^2 < \infty$ for all $\theta, \eta_1$, $\mathbb{E}_{\eta_1}(\|\tilde{U}\|^4) < \infty$, the risks are finite,  $\|g(\tilde{x})\|^2 + 2 \hbox{div} \, g(\tilde{x}) \leq 0$ for all $\tilde{x} \in \mathbb{R}^{d}$, and $0<\alpha< \frac{1}{1+\beta}$.
\end{theorem}
\noindent {\bf Proof.}  We can set $\eta_1=1$ without loss of generality.
The difference in risks is given by
\begin{eqnarray*}
\Delta &=&  \mathbb{E} \left[ \rho (\|\tilde{X} +  \frac{\alpha \|\tilde{U}\|^2}{k+2} \, g(\tilde{X})-\tilde{Y}\|^2 + \|\tilde{U}\|^2) \, - \, \rho (\|\tilde{X} - \tilde{Y}\|^2 + \|\tilde{U}\|^2) \right]  \\
 \, & \leq & \mathbb{E} \left[ \rho'(\|\tilde{X} - \tilde{Y}\|^2 + \|\tilde{U}\|^2) \{  \frac{\alpha^2 (\|\tilde{U}\|^2)^2}{(k+2)^2} \, \|g(\tilde{X})\|^2 + \frac{2\alpha\|\tilde{U}\|^2}{k+2} \, g(\tilde{X})^T(\tilde{X}-\tilde{Y})   \}  \right],
 \end{eqnarray*}
by virtue of the inequality $\rho(A+b) - \rho(A) \leq \rho'(A) \, b$ for concave $\rho$.  With the change of variables $Z=\tilde{X}-\tilde{Y}$, $W=\tilde{Y}+\beta \tilde{X}$, $\tilde{U}=\tilde{U}$ so that $(Z,W,\tilde{U})$ has joint density proportional to $f(\|z\|^2/(1+\beta) + \|\tilde{u}\|^2/(1+\beta) +  (\|w- (1+\beta)\mu\|^2)/(\beta+\beta^2))$, and conditioning on $W$, we have
\begin{equation}
\nonumber
\Delta \,\leq \, \mathbb{E} \left\lbrace \mathbb{E} \left[  \rho'(\|Z\|^2 + \|\tilde{U}\|^2) \, \{  \frac{\alpha^2 (\|\tilde{U}\|^2)^2}{(k+2)^2} \, \|g(\frac{\|Z+W\|}{1+\beta}\|^2 + \frac{2\alpha \|\tilde{U}\|^2}{k+2} \, g(\frac{Z+W}{1+\beta})^T Z \}  \right]  \vline W\right\rbrace.   
\end{equation}
We proceed by showing that the given conditions imply that the inner conditional expectation, given $W=w$ and  denoted $\Delta(w)$, which is taken with respect to the conditional density $f_w(\|z\|^2 + \|\tilde{u}\|^2) \propto f((\|z\|^2/(1+\beta) + (\|\tilde{u}\|^2/(1+\beta) +  (\|w- (1+\beta)\mu\|^2)/(\beta+\beta^2) $, is non-positive for all $w$.  Applying Lemma \ref{Steingeneral} twice for density $f_w$ and associated $F_w$, we obtain

\begin{eqnarray*}
 \Delta(w) & \propto &\int_{\mathbb{R}^{d+k}}  F_w(\|z\|^2 + \|\tilde{u}\|^2) \{\alpha^2 \, \frac{\hbox{div}(\|\tilde{u}\|^2 \tilde{u}}{(k+2)^2} \,\| g(\frac{z+w}{1+\beta}) \|^2 \, + \,  \frac{2\alpha\|\tilde{u}\|^2}{k+2} \, \hbox{div}_z g(\frac{w+z}{1+\beta}) \} \, d\tilde{u} \, dz \\
 &=  &  \alpha \int_{\mathbb{R}^{d+k}}  F_w(\|z\|^2 + \|\tilde{u}\|^2) \, \frac{\alpha \|\tilde{u}\|^2}{k+2} \{ \alpha \| g(\frac{z+w}{1+\beta}) \|^2 \, + \frac{2}{1+\beta} \,  \hbox{div}_z g(\frac{w+z}{1+\beta})  
\} \, d\tilde{u} \, dz, \\
& \leq & 0\,
\end{eqnarray*}
for all $w \in \mathbb{R}^{d}$, given the given conditions on $g$ and $\alpha$.  This completes the proof. \qed \\

\noindent  The above dominance result is wide ranging and is doubly robust.  First,
the class of dominating predictors is vast.  It includes usual shrinkage estimators which satisfy the familiar differential inequality for minimaxity in a point estimation framework.  These include James-Stein, \blue{James-Stein positive-part}, Baranchik type estimators, Bayesian estimators with respect to a superharmonic prior, etc. 
Secondly, the dominance holds simultaneously for a large collection of losses $\rho$, including squared error penalty $\|\delta-Y'\|^2$, other $L^p$ losses with $\rho(t)=t^p$ and $0<p<1$, the case $\rho(t)=\log(t)$ which will serve for our predictive density estimation framework, and many bounded losses such as reflected normal loss $\rho(t)=1- e^{-t/\alpha}$ with $\alpha>0$.    Thirdly, the dominance holds simultaneously for all model densities $f$ provided the risks are finite and $\mathbb{E}_{\theta, \eta_1} \|g(\tilde{X})\|^2 < \infty$.  This includes the normal case with independently distributed $\tilde{X} \sim N_{d}(\mu, I_{d}/\eta_1)$,  $\tilde{Y} \sim N_{d}(\mu, (\beta/\eta_1) I_{d})$, $\tilde{U} \sim N_k(0, ((1+\beta)/\eta_1) I_k)$, as well as scale mixture of normals with $\eta_1$ random for the above triplet.  We point out that the above result does not necessitate that $\rho$ be positive, and negative values for $\rho$ arise naturally for the connected predictive density estimation problem, which we now address with the help of Theorem \ref{main}.

\begin{theorem}
\label{dominancepde}
	Consider model (\ref{modelXUY}) with $d \geq 3$ and the problem of obtaining a predictive density, based on $(X,U)$, of the conditional density of $Y$ given $(X,U)$.  Consider the Bayes predictive density $\hat{q}_{\pi_0}(\cdot;(X,U)) \sim T_d(k,cX, \sqrt{\frac{(1+c^2) \|U\|^2}{k}})\,,$  competing predictive density estimators $q_{\pi_0, \hat{\theta}}(\cdot;(X,U)) \sim T_d(k,c\hat{\theta}(X,U), \sqrt{\frac{(1+c^2) \|U\|^2}{k}})\,,$ and their efficiency as measured by Kullback-Leibler risk.  Then $q_{\pi_0, \hat{\theta}}(\cdot;(X,U))$ dominates $\hat{q}_{\pi_0}(\cdot;(X,U))$ with  $\hat{\theta}(X,U) = X + a \, \frac{ \|U\|^2}{k+2} \, g(\frac{X}{c})$, provided $\mathbb{E}_{\theta, \eta_1} \|g(X)\|^2 < \infty$, $\mathbb{E}(\|U\|^4) < \infty$, finiteness of risk,  $\|g(t)\|^2 + 2 \hbox{div} \, g(t) \leq 0$ for all $t \in \mathbb{R}^{d}$, and $0<a< \frac{(1+c^2)}{c^2 (1+c)}$.
\end{theorem}
\noindent {\bf Proof.}  We make use of Lemma \ref{duality} and Theorem \ref{main}.  With Lemma \ref{duality}'s duality result,  $q_{\pi_0, \hat{\theta}}$ will dominate $\hat{q}_{\pi_0}$ if and only if $c \, \hat{\theta}(X,U)$ dominates $cX$ as a predictor of $Y$ under prediction loss  $\rho_0(\|Y-c\hat{\theta}\|^2 + (1+c^2) \|U\|^2)$   with $\rho_0(t)=\log(t)$ and $c \hat{\theta}$ a given prediction.   With the change of variables
$$   (X,Y,U) \to (\tilde{X}=\frac{X}{c}, \tilde{Y}=\frac{Y}{c^2}, \tilde{U}= \frac{\sqrt{1+c^2}}{c^2} U  )\,, $$
dominance will be achieved if and only if $c \hat{\theta}\left(c\tilde{X}, \frac{c^2}{\sqrt{1+c^2} } \tilde{U}  \right)$ dominates $c^2 \tilde{X}$ as a predictor of $c^2 \tilde{Y}$ under loss
\begin{eqnarray}
\nonumber
\, & \,\;\; & \rho_0 \left(  \| c^2 \tilde{Y} - c^2 \frac{\hat{\theta}(c\tilde{X}, \frac{c^2}{\sqrt{1+c^2}} \tilde{U})}{c} \|^2 + \|c^2 \tilde{U}\|^2 \right)   \\
\, & = & 
\label{newloss} \rho \left(  \| \tilde{Y} - \delta(\tilde{X},\tilde{U}) \|^2 + \| \tilde{U}\|^2 \right) \,,
\end{eqnarray}
with $\rho(t)=\rho_0(c^4 t) \,(= 4 \log c + \log t)\,$ for $t>0,$ and 
$$ \delta(\tilde{X},\tilde{U}) = \frac{\hat{\theta}(c\tilde{X}, \frac{c^2}{\sqrt{1+c^2}} \, \tilde{U})}{c}\,. $$
Dominance will thus be achieved if the above $\delta(\tilde{X},\tilde{U})$ dominates $\tilde{X}$ as a predictor of $\tilde{Y}$ under loss (\ref{newloss}).  Since the triplet $(\tilde{X},\tilde{Y},\tilde{U})$ has density as in (\ref{modelX'Y'U'}) with $\mu= \theta/c$, $\beta=1/c$, $\eta_1 = c^2 \eta$, we can apply Theorem \ref{main}.  Hence, if $\hat{\theta}(X,U)= X + a \, \frac{ \|U\|^2}{k+2} g(\frac{X}{c})$, we have corresponding $\delta(\tilde{X},\tilde{U}) = \tilde{X} + a  \frac{c^3}{1+c^2}\frac{ \|\tilde{U}\|^2}{k+2} g(\tilde{X})$, and a sufficient condition for dominance is indeed
$$ 0 <  a \frac{c^3}{1+c^2} < \frac{1}{1+\beta} \Longleftrightarrow  0 < a < \frac{(1+c^2)}{c^2 (1+c)}\,, $$
since $\beta=1/c$.   \qed \\

	We conclude this section by pointing out that above dominance holds simultaneously for all $f$ subject to the finiteness conditions.

\section{Bayesian representations and robustness results}

\subsection{On posterior robustness under separable priors}

We expand here on a general robustness property where a class of Bayesian inference methods are robust with respect to a model density.  Consider the following canonical set-up represented by spherically symmetric densities, with residual vector $U$,

\begin{equation}
\label{modelXU}
X,U|\theta, \eta \sim \eta^{(d+k)/2} \, f\left(\eta(\|x-\theta|^2+ \|u\|^2 )  \right)\,,
\end{equation}
with $x,\theta \in \mathbb{R}^{d}$, $ u \in \mathbb{R}^{k}$,  $\eta^{-1/2}$ a scale parameter, and $f(\|t\|^2)$ a spherically symmetric density on $\mathbb{R}^{d+k}$. 
Further consider Bayesian inference for separable priors of the form:

\begin{equation}
\label{prior2}  \theta,\eta \sim \pi_1(\theta) \, \eta^a\;; \theta \in \mathbb{R}^{d}, \eta>0, a \in \mathbb{R};
\end{equation}
with $\pi_1(\theta)$ absolutely continuous with respect to a $\sigma-$finite measure $\nu$.  We point out that these priors are necessarily improper.  Whenever the posterior distribution of $(\theta, \eta)$ is well-defined, we have the following general representation.  

 \begin{theorem}
\label{independence}
Consider model (\ref{modelXU}), a prior distribution as in (\ref{prior2}) and, 
for a given $(x,u)$, $\tau= \eta \left(\|\theta-x\|^2 + \|u\|^2  \right)$.  
Assume that  
\begin{equation}
\nonumber
\int_{\mathbb{R}_+} t^{a+ \frac{d+k}{2}} \, f(t) \, dt < \infty \, \hbox{ and } 
\int_{\mathbb{R}^d}  \frac{\pi_1(z)}{\left(\|z-x\|^2 + \|u\|^2  \right)^{a+1 + \frac{d+k}{2}}} \, dz < \infty\,.
\end{equation}
Then, conditional on $(x,u)$, $\theta$ and $\tau$ are independent with 
densities 
\begin{equation}
\label{jointposterior}
\tau|x,u \propto \tau^{a+ \frac{d+k}{2}} f(\tau)   \, \hbox{ and } 
\theta|x,u \propto 
\frac{\pi_1(\theta)}{\left(\|\theta-x\|^2 + \|u\|^2  \right)^{a+1 + \frac{d+k}{2}}}\,\,. 
\end{equation}
Moreover, the marginal posterior distribution of $\theta$ is independent of $f$, while the marginal posterior distribution of $\tau$ is independent of $\pi_1$.
\end{theorem} 
\noindent {\bf Proof.}  We have, for the given model and prior, the posterior density 
\begin{equation}
\pi_{1,a}(\theta, \eta|x,u) \propto \eta^{a+ \frac{d+k}{2}} \, f\left(\eta(\|x-\theta|^2+ \|u\|^2 )  \right) \, \pi_1(\theta)\,.
\end{equation}
The change of variables $(\theta, \eta) \to (\theta, \tau)$ yields (\ref{jointposterior}), with the densities well defined given the finiteness assumption. 
Finally, the conditional independence and posterior marginal distributions follow from (\ref{jointposterior}).   \qed \\

\noindent 
The above independence representation now paves the way to the following results.

\begin{corollary}
\label{independencecorollary}
Consider model (\ref{modelXU}) and a prior distribution as in (\ref{prior2}) for which the posterior distribution of $(\theta, \eta)$ is well-defined.  Then,

\begin{enumerate}
\item[ {\bf (a)}]   Bayesian posterior inference about $\theta$, based solely on the posterior distribution of $\theta$, such as Bayesian confidence regions, tests, and predictors, as well as Bayes point estimators such as $\mathbb{E}(\theta|X,U)$, do not depend on the model density $f$;

\item[ {\bf (b)}]   For $\int_{\mathbb{R}_+} t^{a+b + \frac{d+k}{2}} \, f(t) \, dt < \infty \,$,   Bayes point estimators of $\theta$ under losses of the form  $\eta^b \, \rho(\|\delta-\theta  \|^2)$ are, provided they exist, independent of the model density $f$;

\item[ \bf (c)]   In particular for $b=1, \rho(t)=t$, corresponding to scale invariant squared-error loss $\eta \, \|\delta - \theta\|^2$, the Bayes point estimator of $\theta$ is given by:

\begin{equation}
\label{scaleinvariant}
\delta_{\pi_1,a}(x,u) \, = \,  \frac{\int_{\mathbb{R}^d}  \frac{\theta}{\left((\|\theta - x\|^2 + \|u\|^2 )^{a+2+ \frac{d+k}{2}}  \right)} \;  \pi_1(\theta) \, d\nu(\theta) }{\int_{\mathbb{R}^d}  \frac{1}{\left((\|\theta - x\|^2 + \|u\|^2 )^{a+2+ \frac{d+k}{2}}  \right)} \;  \pi_1(\theta) \, d\nu(\theta)}  \, ,  
\end{equation}
provided that 
$\int_{\mathbb{R}_+} t^{a+1+ \frac{d+k}{2}} \, f(t) \, dt < \infty $ and that 
$\mathbb{E} \!\left( \frac{\|\theta\|^{\ell}}{\|\theta-x\|^2 + \|u\|^2 } \, | \, x,u \right) < \infty$. 
\end{enumerate}
\end{corollary}
\noindent {\bf Proof.}    Part {\bf (a)} follows immediately from Theorem \ref{independence}.  For part {\bf (b)}, the expected posterior loss associated with point estimate $\delta$ is given by (recall 
$\tau= \eta \left(\|\theta-x\|^2 + \|u\|^2  \right)$)
\begin{eqnarray*}
\mathbb{E} \left( \eta^b \, \rho(\|\delta-\theta  \|^2) \, |\, x,u \, \right) & = & \mathbb{E} \left( \tau^b \,  \frac{\rho(\|\delta-\theta  \|^2)}{(\|\theta - x\|^2 + \|u\|^2 )^b} \, | \, x,u \, \right) \\
\, & = & \mathbb{E} \left( \tau^b \, |x,u \,\right) \,\; \mathbb{E} \left( \frac{\rho(\|\delta-\theta  \|^2)}{(\|\theta - x\|^2 + \|u\|^2 )^b} \, |\, x,u \, \right) \, ,
\end{eqnarray*}
with the given finiteness assumption and by making use of Theorem \ref{independence}.  It is thus the case that the minimizing $\delta_{\pi_1,a}(x,u)$ depends only on the posterior distribution of $\theta|x,u$, and consequently is independent of $f$.  Finally, for part {\bf (c)}, we have:

\begin{eqnarray*}
\delta_{\pi_1,a}(x,u) & = &  argmin_{\delta} \,\; \mathbb{E} \left( \frac{\|\delta-\theta\|^2}{\|\theta-x\|^2 + \|u\|^2} \, | x,u  \right) \\
& = &  \frac{\mathbb{E} \left( \frac{\theta}{\|\theta-x\|^2 + \|u\|^2 } \, | \, x,u \right)  } {\mathbb{E} \left( \frac{1}{\|\theta-x\|^2 + \|u\|^2 } \, | \, x,u \right)} \,,
\end{eqnarray*}
by a familiar weighted squared-error loss Bayes estimator representation.   The result then follows by incorporating the posterior density given in (\ref{jointposterior}).  \qed \\

\noindent   The above results are indeed quite striking and analogous results for predictive densities will be elaborated on below.   However, from a historical perspective, the findings above add to, extend, or clarify earlier findings.  More precisely, the robustness of the point estimators in (\ref{scaleinvariant}) was observed by Maruyama (2003) for $\pi_1(\theta)$ of the form $\|\theta\|^{b}$,  Fourdrinier and Strawderman (2010) as well as Maruyama and Strawderman (2005) for further separable priors, and Jafari Jozani, Marchand and Strawderman (2013) for the univariate case of a positive $\theta$ with $\pi_1(\theta)=\mathbb{I}_{(0,\infty)}(\theta)$. The results of Theorem \ref{independence} and Corollary \ref{independencecorollary} are much more general though.  This includes numerous possible forms of $\pi_1$.  For instance, in the univariate case $d=1$, and with the restriction $\theta \in [-m,m]$, and the two-point uniform boundary prior  (i.e., $\pi_1(m)=\pi_1(-m)=1/2$), expression (\ref{scaleinvariant}) yields  $m (B-A)/(B+A) $ with $B= \left\lbrace (x+m)^2 + u^2 \right\rbrace^{a+(k+5)/2}$ and  $A= \{(x-m)^2 + u^2\}^{a+(k+5)/2} $. \\

\noindent  Although the focus of this paper is not on frequentist risk comparisons of point estimators, we conclude with a robust dominance result illustrating how naturally a dominance finding in the normal case can carry-over to dominance findings for scale mixture of normals.

\begin{theorem}
Consider model (\ref{modelXU}) and the problem of estimating $\theta$ based on $(X,U)$ and loss $L((\theta, \eta), \hat{\theta})$.  Suppose that $\hat{\theta}_1(X,U)$ dominates 
$\hat{\theta}_0(X,U)$ with smaller expected loss for all $(\theta, \eta)$.  Then, $\hat{\theta}_1(X,U)$ also dominates $\hat{\theta}_0(X,U)$ for all scale mixture of normals $f$ as long as the corresponding risks are finite.   
\end{theorem}
\noindent  {\bf Proof.}   We have the representation  $(X,U)|Z \sim N_{d}(\theta, (Z/\eta) I_{d})$ that permits to write the difference in risks as equal to
$$   \mathbb{E}^Z \left\lbrace  \, \mathbb{E}^{X,U|Z}  \left( L((\theta, \eta), \hat{\theta}_1(X,U)) - L((\theta, \eta), \hat{\theta}_0(X,U)) \right)   \right\rbrace\,. $$
The result follows since the inner expectation is negative with probability one (with respect to $Z$) by virtue of the dominance assumption in the normal case.  \qed

\subsection{On a predictive density estimation representation and robustness property}

We follow-up with a further robustness result applicable in the predictive density estimation framework presented in part {\bf C.} of the Introduction. 

Reconsider model (\ref{modelXUY}) and the problem of obtaining a predictive density for
the conditional density $q_{\theta, \eta}(y|x,u), y \in \mathbb{R}^{d}$ and assessing its efficiency with respect to Kullback-Leibler loss (\ref{Kullback-Leibler}).  Now, for separable priors as in (\ref{prior2}), we have the following robustness property.

\begin{theorem}  
\label{independentoff}
Suppose $\int_{\mathbb{R}_+} t^{d+k/2+a}  \,  f(t)  \, dt  \, < \infty$. 
Assuming the posterior distribution of $\theta, \eta$ is well-defined, Bayesian predictive densities under Kullback-Leibler loss, for prior densities of the form $\theta,\eta \sim \pi_1(\theta) \, \eta^a\,$, are independent of the model density $f$ and given by
\begin{equation}
\label{pdeformula}
\hat{q}_{\pi_1}(y;x,u) \,=\, \frac{\int_{\mathbb{R}^{d}} (\|x-\theta \|^2 + \|u\|^2  + \|y-c\theta \|^2 )^{-n} \, \pi_1(\theta) \, d\nu(\theta) }  {\int_{\mathbb{R}^{2d}}  (\|x-\theta \|^2 + \|u\|^2  + \|y-c\theta \|^2 )^{-n} \, \pi_1(\theta) \, d\nu(\theta) \, dy }  \,
\end{equation}
with $n=d+k/2+a+1\,$.
\end{theorem}
\noindent  {\bf Proof.}  
From Aitchison (1975), we have
\begin{equation}
\hat{q}_{\pi_1}(y;x,u) \,= \, \int_{\mathbb{R}^{d+1}} q_{\theta, \eta}(y|x,u) \, 
\pi(\theta, \eta|x,u) \, d\nu(\theta) \, d\eta \,.
\end{equation}
As in Theorem \ref{independence}, we re-express the posterior in terms of $\theta, \tau$, with $\tau= \eta (\|x-\theta \|^2 + \|u\|^2  + \|y-c\theta \|^2 )$, to obtain
\begin{equation}
\hat{q}_{\pi_1}(y;x,u) \, \propto \, \int_{\mathbb{R}^{d+1}} \tau^{d+a+k/2} f(\tau) \,\,
\frac{\pi_1(\theta)}{(\|x-\theta \|^2 + \|u\|^2  + \|y-c\theta \|^2)^n} \, d\tau \, d\nu(\theta) \,.
\end{equation}
This now yields the result as the term $\int_{\mathbb{R}_+} \tau^{d+a+k/2} f(\tau) \, d\tau $
is constant and factors in both the numerator and denominator of (\ref{pdeformula}).   \qed \\

\noindent The minimum risk equivariant predictive density $\hat{q}_{\pi_0}$ solution, previously stated in (\ref{mre}), is obtained as a particular case of (\ref{pdeformula}) with $\pi_1(\theta)=1$, $\nu$ the Lebesgue measure on $\mathbb{R}^{d}$, and $a=-1$.  It can be computed directly, or inferred from the normal case solution (e.g.., Aitchison and Dunsmore, 1975; Kato, 2009).   Illustrating the former,  as a particular case of priors
(\ref{prior2}) with $\pi_1 \equiv 1$, we have setting $B= ( \frac{\| y - cx \|^2}{(1+c^2) \|u\|^2}+ 1 )$ and with the decomposition  $\|x-\theta\|^2 + \|y-cx\|^2 \, = \, (1+c^2) \, (\|\theta - (\frac{x+cy}{1+c^2}) \|^2 + \frac{\|y-cx\|^2}{(1+c^2)^2} )$:
\begin{eqnarray*}
\hat{q}_{\pi_0,a}(y;x,u) & \propto & \int_{\mathbb{R}^{d}}  (\|x-\theta \|^2 + \|u\|^2  + \|y-c\theta \|^2 )^{-n}  \, d\theta \\
\, & \propto &   B^{- (d+k+2a+2)/2}\int_{\mathbb{R}^{d}} B^{-d/2} \,  \left( \frac{\|\theta - (\frac{x+cy}{1+c^2}) \|^2}{B} + 1 \right)^{-(d/2 + (d+k+2a+2)/2)} \, d\theta \,\\
\, & \propto &  ( \frac{\| y - cx \|^2}{(1+c^2) \|u\|^2}+ 1 )^{- (d+k+2a+2)/2}\,,
\end{eqnarray*}
which is a Student $T_d(k+2a+2, cx, \sqrt{\frac{(1+c^2) \|u\|^2}{k+2a+2}})$ density, and which yields (\ref{mre})  indeed for $a=-1$. 
 
\section{A Bayesian dominance result for scale mixture of normals}

We consider here Bayes predictive densities $\hat{q}_{\pi_1}$ with respect to separable prior densities $\theta, \eta \sim \pi_1(\theta) \, \eta^{-1}$ and comparisons with the particular case $\hat{q}_{mre}=\hat{q}_{\pi_0}$, which is a Bayes predictive density for prior density $\theta, \eta \sim \eta^{-1}$.  The next result, stated a little more generally, will imply that a dominance result which holds in the normal case $f(t)= \phi(t)$ in (\ref{modelXUY}) will necessarily hold simultaneously for all variance  mixture of normals densities with $f(t) = \int_{\mathbb{R}_+} z^{-d} \phi(z^{-1}t) \, dG(z)$, $G$ being the c.d.f. of the mixing scale parameter.

\begin{theorem}
\label{normalimplies...}
Consider model (\ref{modelXUY}) and the problem of obtaining a predictive density for Kullback-Leibler loss (\ref{Kullback-Leibler}), based on $(X,U)$, for the conditional density of $Y$ given $(X,U)$.   Then, subject to the finiteness of risks, if $\hat{q}_{1}$ dominates  $\hat{q}_{0}$ in the normal case, then $\hat{q}_{1}$ dominates $\hat{q}_{0}$ simultaneously for all scale mixtures of normals.
\end{theorem}
\noindent {\bf Proof.}    
We have from (\ref{riskKL}), for the difference in risks,
\begin{eqnarray*}
R_{KL}((\theta, \eta), \hat{q}_{0}) & - & R_{KL}((\theta, \eta), \hat{q}_{1}) \, \\ \, & = & \mathbb{E}^Z   \mathbb{E}^{X,U,Y|Z}  \log \left(\frac{ \hat{q}_1(Y;X,U)}{\hat{q}_{0}(Y;X,U)}   \right) \,  \\
\, & = & \mathbb{E}^Z  \Delta(\theta, \eta, Z) \; \;  (\hbox{ say }),
\end{eqnarray*}
with $X,U,Y|Z$ normally distributed, as in (\ref{modelXUY}) with $f(t) = z^{-(2d+ k)} \phi(z^{-1} t )$, and $Z$ having c.d.f. $G$ on $(0,\infty)$.   Now, the assumptions imply that $\Delta(\theta, \eta, Z) \geq 0$ for all $\theta \in \mathbb{R}^{d}, \eta>0$ with probability one, with strict inequality for some $(\theta, \eta)$, thus establishing the result.    \qed

\begin{corollary}
\label{Katomixtures}
	Consider model (\ref{modelXUY}) with a scale of mixtures of normals $f$, $d \geq 3$, and the problem of obtaining a predictive density, based on $(X,U)$, of the conditional density of $Y$ given $(X,U)$.  Then, the Bayes predictive density estimator $\hat{q}_{\pi_h}$ with respect to the harmonic prior density $\pi_h(\theta, \eta) = \eta^{-1} \|\theta\|^{2-d}$ dominates the MRE predictive density $\hat{q}_{\pi_0}$ under Kullback-Leibler loss.  Furthermore, the dominance holds simultaneously for all scale mixture of normals $f$.
\end{corollary}
\noindent {\bf Proof.} 
With $\hat{q}_{\pi_h}$ dominating $\hat{q}_{\pi_0}$ in the normal case by virtue of Kato (2009), since both $\hat{q}_{\pi_h}$ and $\hat{q}_{\pi_0}$ do not vary with $f$ by Theorem \ref{independentoff}, the result is a direct consequence of Theorem \ref{normalimplies...}.  \qed

\section*{Acknowledgements}

Dominique Fourdrinier's research is partially supported by the RNF (Russian 
National Foundation), Project \#17-11-01049 and by the Tassili Program, 
Project \#18MDU105.  \'Eric Marchand's research is supported in part by the 
Natural Sciences and Engineering Research Council of Canada, and William 
Strawderman's research is partially supported by grants from the Simons 
Foundation (\#209035 and \#418098).

\renewcommand{\baselinestretch}{0.8}

\end{document}